\documentclass[11pt]{amsart}
\usepackage{amsmath,amssymb,pxfonts}
\title{Discretization Error of Stochastic Integrals}
\author{Masaaki Fukasawa}
\thanks{
Masaaki Fukasawa\\
Center for the Study of Finance and Insurance, Osaka University\\
Japan Science and Technology Agency, CREST\\ 
1-3 Machikaneyama, Toyonaka, Osaka, Japan}
\date{2009-12-21}

\newtheorem{thm}{Theorem}[section]
\newtheorem{prop}[thm]{Proposition}
\newtheorem{lem}[thm]{Lemma}

\newtheorem{defn}[thm]{Definition}
\newtheorem{assu}[thm]{Condition}

\begin{document}
\maketitle
\begin{abstract}
Asymptotic error distribution for approximation
of a stochastic integral with respect to continuous semimartingale
by Riemann sum with general stochastic
partition is studied.
Effective discretization schemes of which
asymptotic conditional mean-squared error
 attains a lower bound are constructed.
Two applications are given;
efficient delta hedging strategies with transaction costs
and effective discretization schemes for the Euler-Maruyama
approximation are constructed.
\end{abstract}




\section{Introduction}
The present article  studies the asymptotic distribution of
a sequence of continuous processes $Z^n = \{Z^n_t\}_{t \in [0,T)}$ defined as
\begin{equation} \label{defz}
Z^n_t = \int_0^t X_s\mathrm{d}Y_s - \sum_{j=0}^{\infty}
X_{\tau^n_j}(Y_{\tau^n_{j+1} \wedge t}-Y_{\tau^n_j \wedge t})
\end{equation}
for one-dimensional continuous semimartingales $X = \{X_t, \mathcal{F}_t\}$, 
$Y=\{Y_t, \mathcal{F}_t\}$ and
sequences of $\{\mathcal{F}_t\}$-stopping times $\tau^n = \{\tau^n_j\}$ 
with
\begin{equation} \label{tau}
0 = \tau^n_0 < \tau^n_1 < \dots < \tau^n_j < \dots, \ \ 
\lim_{j\to \infty} \tau^n_j = T \text{ a.s.,}
\end{equation}
where $T \in (0,\infty]$ is fixed and
we suppose that the intervals $\tau^n_{j+1}-\tau^n_j$ 
converge to $0$ as $n \to \infty $ in a sense specified later.
Taking into mind the definition of stochastic integrals, the 
asymptotic behavior of $Z^n$ is theoretically of interest.
There may be no need to explain its practical importance.
This problem was studied by Rootz\'en~\cite{Rootzen} in
the case that $Y$ is a Brownian motion and
the asymptotic distribution was specified in the case that
$\tau^n_j = j/n$ and $X_s = f(Y_s,s)$ for a smooth function $f$.
Jacod~\cite{Jacod1994} treated a related problem
on the condition that each interval $\tau^n_{j+1} - \tau^n_j$ is
$\mathcal{F}_{\tau^n_j}$ measurable.
Jacod and Protter~\cite{JP} considered the case $X=Y$ and 
$\tau^n_j = j/n$ and derived the asymptotic error distribution of
the Euler-Maruyama scheme for stochastic differential equations.
Hayashi and Mykland~\cite{HM2005} discussed this problem
again for the case $\tau^n_j = j/n$ in the context of
discrete-time hedging error.
Geiss and Toivola~\cite{GT} treated an irregular deterministic
discretization scheme.
The condition $\tau^n_j = j/n$, or more generally,
that each interval  $\tau^n_{j+1} - \tau^n_j$ is
$\mathcal{F}_{\tau^n_j}$ measurable, played an important role in
those preceding studies of central limit theorem.
On the other hand, 
Karandikar~\cite{Karandikar1995}
constructed a discretization scheme $\tau^n$
such that $Z^n_t$ converges to $0$ almost surely.
Since the almost sure convergence is not usually attained by the
time-equidistant scheme  $\tau^n_j = j/n$, 
it means Karandikar's scheme is more effective in a sense.
The scheme is defined using passage times of $X$, so that
the measurability condition, that each interval $\tau^n_{j+1} - \tau^n_j$ is
$\mathcal{F}_{\tau^n_j}$ measurable, is not satisfied.
Recently, Fukasawa~\cite{Fukasawa2009}
gave the asymptotic distribution for such a scheme and
Fukasawa~\cite{F2009} extended the results
to a more general scheme in the case $X=Y$.
The present article extends those limit theorems  
and constructs effective discretization schemes
of which asymptotic conditional mean-squared error 
 attains a lower bound. In particular, Karandikar's scheme is 
shown to be superior to the time-equidistant one also in
terms of mean-squared error.
An application to delta hedging with transaction costs is given,
which can be directly used in practice.
It remains for further research to extend the results
to discontinuous semimartingales.
In Section~2, we describe and prove main results. 
Effective discretization schemes are constructed in Section~3.
The application to hedging is presented in Section~4.
An application to the Euler-Maruyama approximation is given in Section~5.

\section{Central Limit Theorem}
\subsection{Notation and conditions}
Here we give a rigorous formulation and 
describe several conditions on $X$, $Y$ and $\tau^n$.
Let $(\Omega, \mathcal{F}, \{\mathcal{F}_t\}_{t \geq 0}, P)$ be a 
filtered probability space.
The filtration $\{\mathcal{F}_t\}$ is assumed to satisfy 
the usual conditions.
We denote by $F_1 \cdot F_2$ the Stieltjes integral or the stochastic
integral of $F_1$ with respect to $F_2$.

Let us recall the definition of stable convergence.
Let $\mathbb{E}$ be a complete separable metric space and 
$F^n$ be a sequence of $\mathbb{E}$-valued random variables
defined on $(\Omega, \mathcal{F}, P)$.
\begin{defn}
For a sub $\sigma$-field  $\mathcal{G} \subset \mathcal{F}$,
we say $F^n$ converges $\mathcal{G}$-stably if
for all $\mathcal{G}$-measurable random variable $F_0$,
the joint distribution  $(F^n, F_0)$ converges in law.
\end{defn}

Our main results stated in the next subsection
are stable convergences of $Z^n$ defined as (\ref{defz}) with
continuous semimartingales $X$ and $Y$.
Notice that a stable convergence is stable against, in particular,
the usual localization procedure as well as 
the Girsanov-Maruyama transformation.

Denote by $\mathcal{P}$ and $\mathcal{P}_0$ 
the set of the predictable processes and the set of
locally bounded left-continuous adapted processes respectively.
Let $T \in (0,\infty]$ be fixed. 
Given a continuous semimartingale $M$ and $k\in \mathbb{N}$, put
\begin{equation*}
\mathcal{P}_M^k = \left\{ H \in \mathcal{P}; 
|H|^k \cdot \langle M \rangle_t < \infty \text{ for all } t \in [0,T)
\right\}.
\end{equation*}
Denote by $\mathcal{S}$ the set of
the continuous semimartingales $(X,Y,M)$ satisfying the
following Condition~\ref{YN}.
\begin{assu} \label{YN}
There exist 
$\psi, \varphi, \kappa \in \mathcal{P}_M^2$
and a locally bounded predictable process $\gamma \in \mathcal{P}$ 
such that
\begin{equation*}
X  =X_0 +  \psi \cdot \langle M \rangle + \gamma \cdot M, \ \ 
Y = Y_0 + \varphi \cdot \langle M \rangle + M^Y
\end{equation*}
on $[0,T)$,
where $M^Y$ is a  continuous local martingale with
\begin{equation*}
\langle M^Y \rangle = \kappa \cdot \langle M \rangle.
\end{equation*}
In addition, $M$ is a continuous local martingale with  
$E[\langle M \rangle_T^6] < \infty$.
\end{assu}
The integrability of $\langle M \rangle_T$ is not restrictive in light of
the localization procedure.
In order to describe 
conditions on $\tau^n$, put
\begin{equation*}
G_{j,n}^1 = E[|M_{\tau^n_{j+1}}- M_{\tau^n_j}|| 
\mathcal{F}_{\tau^n_j}], \ \
G_{j,n}^k = E[(M_{\tau^n_{j+1}}- M_{\tau^n_j})^k| 
\mathcal{F}_{\tau^n_j}] 
\end{equation*}
for a given continuous local martingale $M$ with
$E[\langle M \rangle_T^6] < \infty$ and
$k \in \mathbb{N}$ with $2 \leq k \leq 12$.
In addition, put
\begin{equation} \label{ntau}
N[\tau^n]_{\tau} = \max\{j \geq 0; \tau^n_j \leq \tau\}
\end{equation}
for a given stopping time $\tau$. 
Denote by $\mathcal{T}(M)$ the set of 
the sequences of stopping times $\{\tau^n\}$ satisfying
(\ref{tau}) and the following Condition~\ref{gs}.
\begin{assu} \label{gs}
There exist a sequence $\epsilon_n$ with $\epsilon_n \to 0$
and $a, b \in \mathcal{P}_0$ such that
\begin{equation*}
G_{j,n}^4/G_{j,n}^2 = a_{\tau^n_j}^2 \epsilon_n^2
+ o_p(\epsilon_n^2), \ \ 
G_{j,n}^3/G_{j,n}^2 = b_{\tau^n_j}\epsilon_n
+ o_p(\epsilon_n)
\end{equation*}
and
\begin{equation*}
G_{j,n}^6/G_{j,n}^2 = O_p(\epsilon_n^{4}), \ \ 
G_{j,n}^{12}/G_{j,n}^2 = o_p(\epsilon_n^{8}), \ \ 
\end{equation*}
uniformly in $j=0,1,\dots, N[\tau^n]_t$
for all $t \in [0,T)$.
\end{assu}
Condition~\ref{gs} is 
slightly stronger than  Condition~1 of Fukasawa~\cite{F2009};
nevertheless, all examples given by Fukasawa~\cite{F2009}
satisfy also this condition.
Here, $\epsilon_n$ serves as a scale of increments of $M$.
Note that $G_{j,n}^4/G_{j,n}^2 = O_p(\epsilon_n^2)$
implies $G_{j,n}^3/G_{j,n}^2 = O_p(\epsilon_n)$
by Lemma~\ref{PearsonLem}.
In usual cases, we have $G_{j,n}^{2k} = O_p(\epsilon_n^{2k})$,
which in fact holds, for example, if 
$\mathrm{d}\langle M \rangle_t/\mathrm{d}t$ exists
and is bounded and if $\tau^n_{j+1}-\tau^n_j$ is of
$O(\epsilon_n^2)$ uniformly.
Condition~\ref{gs} is therefore a quite mild condition
in the context of high-frequency asymptotics.
It is often easily verified 
by using the Dambis-Dubins-Schwarz 
time-change technique for martingales
when $\tau^n$ is a function of the path of  $M$.
Once it reduces to the Brownian motion case
by the time-change,
one can utilize various results on 
Brownian stopping times.
See Fukasawa~\cite{F2009} for examples.
In light of the Skorokhod stopping problem, 
the distribution of an increment can be any centered distribution
with a suitable moment condition.
The left continuity of $a^2$ and $b$
corresponds to  a local homogeneity property of the distributions of increments.
It should be noted that $\sup_j|\tau^n_{j+1} \wedge t-\tau^n_j \wedge t| \to 0$
does not follow from Condition~\ref{gs} nor needed for
our main results. See Lemma~\ref{sup} for 
what follows instead.

Denote by $\mathcal{T}_1(M)$, $\mathcal{T}_2(M)$ the subsets of 
$\mathcal{T}(M)$
satisfying the following Condition~\ref{gz} and
Condition~\ref{g2} respectively.
\begin{assu} \label{gz}
In addition to Condition~\ref{gs}, there exists
$\zeta \in \mathcal{P}_0$ 
such that 
\begin{equation*}
\zeta^{-1} \in \mathcal{P}_0, \ \ 
\epsilon_n G_{j,n}^1/G_{j,n}^2
 = \zeta_{\tau^n_j} + o_p(1)
\end{equation*}
uniformly in $j=0,1,\dots, N[\tau^n]_t$ 
for all $t \in [0,T)$.
\end{assu}
\begin{assu}\label{g2}
In addition to Condition~\ref{gs}, there exists  $q \in \mathcal{P}_0$ 
such that 
\begin{equation*}
q^{-1} \in \mathcal{P}_0, \ \ 
G_{j,n}^2 = q^2_{\tau^n_j} \epsilon_n^2 + o_p(\epsilon_n^2)
\end{equation*}
uniformly in $j=0,1,\dots, N[\tau^n]_t$ for all $t \in [0,T)$.
\end{assu}

Finally, for $\tau^n$ with (\ref{tau})
and $t \in [0,T)$, put
\begin{equation}\label{kak}
[M]^t_{j,n} = \langle M \rangle_{\tau^n_{j+1} \wedge t}
-  \langle M \rangle_{\tau^n_{j} \wedge t}.
\end{equation}

\subsection{Main results}
Here we state general results on the asymptotic 
distribution of $Z^n$. The proofs are deferred to Section~\ref{proofs}.

\begin{thm} \label{thm1}
Let $(X,Y,M) \in \mathcal{S}$, $\tau^n \in \mathcal{T}(M)$
and $Z^n$ be defined as (\ref{defz}).
Assume one of the following two conditions to hold:
\begin{itemize}
\item $M$ is the local martingale part of $X$, that is, 
\begin{equation}\label{gam1}
\gamma \equiv 1.
\end{equation}
\item for all $t \in [0,T)$,
\begin{equation} \label{add}
E\left[
\sum_{j=0}^{\infty}|[M]_{j,n}^t|^k
\right] = O(\epsilon_n^{2(k-1)})
\end{equation}
for $k\in \{1,2,3,4,5\}$, 
where $[M]^t_{j,n}$ is defined as (\ref{kak}).
\end{itemize}
Then, 
$Z^n/\epsilon_n$ converges $\mathcal{F}$-stably to
the law  of
\begin{equation} \label{Z}
\frac{1}{3} (b\gamma)\cdot Y  + \frac{1}{\sqrt{6}} (c\gamma)\cdot Y^\prime
\end{equation}
as a $C[0,T)$-valued sequence, where
\begin{equation} \label{cY}
c^2 = a^2 - \frac{2}{3}b^2, \ \ 
Y^\prime = W_{\langle Y \rangle}
\end{equation}
and $W$ is a standard 
Brownian motion which is independent to $\mathcal{F}$.
\end{thm}
Note that the asymptotic distribution (\ref{Z}) is an
$\mathcal{F}$-conditionally Gaussian process, so that
the marginal law is a mixed normal distribution.
The following theorems give lower bounds for the
conditional variance of the mixed normal distribution.
\begin{thm} \label{thm5}
Let $(X,Y,M) \in \mathcal{S}$, $\tau^n \in \mathcal{T}_1(M)$
and $Z^n$ be defined as (\ref{defz}).
Let $u \in \mathcal{P}_0$ and put
\begin{equation*}
U^n_t = \sum_{j=0}^{\infty}
 |u_{\tau^n_j}||M_{\tau^n_{j+1} \wedge t}-M_{\tau^n_j \wedge t}|.
\end{equation*}
Then, it holds
\begin{equation} \label{uconv}
\epsilon_nU^n_t \to U_t :=
(|u|\zeta) \cdot \langle M \rangle_t
\end{equation}
in probability for all $t \in [0,T)$.
Moreover, if (\ref{gam1}) or (\ref{add}) holds, then
$U^n Z^n$ converges $\mathcal{F}$-stably to the law of $UZ$
as a $D[0,T)$-valued sequence, 
where $Z$ is defined as (\ref{Z}).
The asymptotic conditional variance $V_t$
of $ U^n_t Z^n_t $ with $t \in [0,T)$ 
satisfies
\begin{equation} \label{boundu2}
V_t = \frac{1}{6}(c\gamma)^2\cdot \langle Y \rangle_t 
|(|u|\zeta)\cdot \langle M \rangle_t|^2 \geq
\frac{1}{6}
\left|
(|u\gamma|^{2/3}\kappa^{1/3})\cdot \langle M \rangle_t\right|^3
 \text{ a.s..}
\end{equation}
\end{thm}

\begin{thm} \label{thm6}
Let $(X,Y,M) \in \mathcal{S}$, $\tau^n \in \mathcal{T}_2(M)$
and $Z^n$ be defined as (\ref{defz}).
Then, it holds that
\begin{equation} \label{Nconv}
\epsilon_n^2 N[\tau^n]_t \to N_t :=q^{-2}\cdot \langle M \rangle_t
\end{equation}
in probability for all $t \in [0,T)$.
Moreover, if (\ref{gam1}) or (\ref{add}) holds, then
$ \sqrt{N[\tau^n]}Z^n $ converges
$\mathcal{F}$-stably to the law of $\sqrt{N}Z$ as a $D[0,T)$-valued sequence, 
where $Z$ is defined as (\ref{Z}).
The asymptotic conditional variance $V_t$
of $ \sqrt{N[\tau^n]_t}Z^n_t $ with $t \in [0,T)$ 
satisfies
\begin{equation} \label{bound2}
V_t = \frac{1}{6}(c \gamma)^2\cdot \langle Y \rangle_t \ 
q^{-2}\cdot \langle M \rangle_t 
\geq \frac{1}{6}\left\{\left(|\gamma| 
\sqrt{\kappa } \right) \cdot \langle M \rangle_t \right\}^2 \text{ a.s..}
\end{equation}
\end{thm}

Note that the right hand sides of (\ref{boundu2}), (\ref{bound2})
do not depend on the discretization scheme $\tau^n$.
In Section~\ref{effective}, we construct schemes
$\tau^n \in \mathcal{T}_1(M) \cap \mathcal{T}_2(M)$ 
which attain the lower bounds (\ref{boundu2}), (\ref{bound2})
respectively. Its practical meaning is discussed in Sections~4 and 5.

The condition (\ref{add}) will be  easily verified 
especially if
$M$ is a Brownian motion.
Condition~\ref{gs} is then also easily verified 
if, in addition,  $\tau^n$ satisfies the condition
that each interval $\tau^n_{j+1}- \tau^n_j$ is
$\mathcal{F}_{\tau^n_j}$ measurable.
It is therefore not difficult to recover the preceding results
from Theorem~\ref{thm1}.
An irregular scheme treated in Geiss and Toivola~\cite{GT}
is an example.

\subsection{Proof for Theorems} \label{proofs}
Here we give proofs for main results stated in the previous 
subsection.
\begin{lem}\label{unify}
Let $M$ be a continuous local martingale with
$E[\langle M \rangle_T^6] < \infty$ and
$\tau^n \in \mathcal{T}(M)$.
Let $H$, $\gamma$ be a locally bounded cag adapted process
and a locally bounded predictable process respectively.
Put $\bar{M} = \gamma \cdot M$ and define $H^n$, $\bar{M}^n$ as
\begin{equation} \label{bars}
H^n_s = H_{\tau^n_j}, \ \ 
\bar{M}^n_s = \bar{M}_{\tau^n_j},
 \text{ for } j\geq0 \text{ with } s \in [\tau^n_j, \tau^n_{j+1}).
\end{equation}
Assume one of the following two conditions:
\begin{itemize}
\item $\gamma \equiv 1$,
\item (\ref{add}) holds.
\end{itemize}
Then, it holds
\begin{equation} \label{Fconv}
\begin{split}
&\epsilon_n^{-1}((\bar{M}-\bar{M}^n) H^n) \cdot \langle M \rangle_t 
\to \frac{1}{3}(bH\gamma) \cdot \langle M \rangle_t, \\
&\epsilon_n^{-2}((\bar{M}-\bar{M}^n)^2 H^n) \cdot \langle M \rangle_t 
\to \frac{1}{6}(a^2H\gamma^2) \cdot \langle M \rangle_t
\end{split}
\end{equation}
 uniformly in $t$ on compact sets of $[0,T)$
in probability. Moreover,
\begin{equation} \label{Fconv2}
\epsilon_n^{-4}((\bar{M}-\bar{M}^n)^4 H^n) \cdot \langle M \rangle_t 
=O_p(1)
\end{equation}
for all $t \in [0,T)$.
\end{lem}
{\it Proof: }
By the usual localization argument, we can assume $H$, $\gamma$, 
$a$, $b$, $M$
and $\langle M \rangle$  are bounded without loss of generality.
Let us suppose $\gamma\equiv 1$.
Then, for any $l \in \mathbb{Z}$,
\begin{equation*}
\begin{split}
&\epsilon_n^{-l}((\bar{M}-\bar{M}^n)^l H^n) \cdot \langle M \rangle_t \\
&= \epsilon_n^{-l} 
\sum_{j=0}^{\infty} H_{\tau^n_j}
\left\{
\alpha_l (M_{\tau^n_{j+1}\wedge t}  -
M_{\tau^n_{j}\wedge t})^{l+2} +
\beta_j \int_{\tau^n_j \wedge t}^{\tau^n_{j+1}\wedge t}
(M_s - M_{\tau^n_j})^{l+1}\mathrm{d}M_s
\right\}
\end{split}
\end{equation*}
by It$\hat{\text{o}}$'s formula, where $\alpha_l$, $\beta_l$
are constants only depending on $l$ and, in particular,
$\alpha_1 = 1/3$, $\alpha_2 = 1/6$.
Since
\begin{equation*}
\begin{split}
& \epsilon_n^{-2l} 
\sum_{j=0}^{N[\tau^n]_t} H_{\tau^n_j}^2
E\left[\int_{\tau^n_j}^{\tau^n_{j+1}}
(M_s - M_{\tau^n_j})^{2l+2}\mathrm{d}\langle M \rangle_s
\big| \mathcal{F}_{\tau^n_j} \right] \\
&= \epsilon_n^{-2l} \alpha_{2l+2}
\sum_{j=0}^{N[\tau^n]_t} H_{\tau^n_j}^2
G_{j,n}^{2l+4}
\to 0
\end{split}
\end{equation*}
for $l \in {1,2,4}$ by Condition~\ref{gs} and Lemma~\ref{sup},
we have
\begin{equation*}
\epsilon_n^{-l} 
\sum_{j=0}^{\infty} H_{\tau^n_j}
\int_{\tau^n_j \wedge t}^{\tau^n_{j+1}\wedge t}
(M_s - M_{\tau^n_j})^{l+1}\mathrm{d}M_s
\to 0
\end{equation*}
as well as
\begin{equation*}
\epsilon_n^{-l} 
\sum_{j=0}^{\infty} H_{\tau^n_j}
(M_{\tau^n_{j+1}\wedge t}  -
M_{\tau^n_{j}\wedge t})^{l+2}
- \epsilon_n^{-l} 
\sum_{j=0}^{N[\tau^n]_t}H_{\tau^n_j}G_{j,n}^{l+2} \to 0
\end{equation*}
by Lemma~\ref{moreless}. The result then follows from Lemma~\ref{sup}.

Next, let us suppose (\ref{add}) and $\gamma \not\equiv 1$. 
Note that for all $\delta > 0$, there exists
a bounded cag process $\gamma^\delta$ such that
\begin{equation*}
E[|\gamma-\gamma^\delta|^k\cdot \langle M \rangle_t] < \delta
\end{equation*}
for $k=4,6,12$ by Lemma~\ref{cag}.
Notice that for $p, q > 1$ with $1/p + 1/q =1$,
\begin{equation*}
\begin{split}
&E \left[\epsilon_n^{-l}\sum_{j=0}^{\infty}
\int_{\tau^n_j \wedge t}^{\tau^n_{j+1} \wedge t} \left|
\int_{\tau^n_j \wedge t}^{s}
(\gamma_u-\gamma^\delta_u)\mathrm{d}M_u \right|^{l} \mathrm{d}\langle M \rangle_s \right] \\
& \leq E\left[\epsilon_n^{-l}\sum_{j=0}^{\infty}
\sup_{s \in [\tau^n_j \wedge t, \tau^n_{j+1} \wedge t]}
\left|
\int_{\tau^n_j \wedge t}^{s}
(\gamma_u-\gamma^\delta_u)\mathrm{d}M_u \right|^{l} [M]_{j,n}^t \right] \\
& \leq 
C\epsilon_n^{-l}
E\left[ \sum_{j=0}^{\infty} \left|
\int_{\tau^n_j \wedge t}^{\tau^n_{j+1} \wedge t} (\gamma_u-\gamma^\delta_u)^2
 \mathrm{d}\langle M \rangle_u
\right|^{pl/2}
\right]^{1/p} E\left[
\sum_{j=0}^{\infty} \left| [M]_{j,n}^t \right|^q
\right]^{1/q}
\end{split}
\end{equation*}
by the H$\ddot{\text{o}}$lder and the Burkholder-Davis-Gundy inequality.
Furthermore,
\begin{equation*}
\begin{split}
&E\left[ \sum_{j=0}^{\infty} \left|
\int_{\tau^n_j \wedge t}^{\tau^n_{j+1} \wedge t} (\gamma_u-\gamma^\delta_u)^2
 \mathrm{d}\langle M \rangle_u
\right|^{pl/2}
\right]  \\
&\leq E\left[ \sum_{j=0}^{\infty} \left|
\int_{\tau^n_j \wedge t}^{\tau^n_{j+1} \wedge t} 
|\gamma_u-\gamma^\delta_u|^{2pl}
 \mathrm{d}\langle M \rangle_u |[M]_{j,n}^t|^{pl-1}
\right|^{1/2}
\right] \\
&\leq
E\left[\int_0^t |\gamma_u-\gamma^\delta_u|^{2pl}  \mathrm{d}\langle M \rangle_u
\right]^{1/2}
E\left[
\sum_{j=0}^{\infty} |[M]_{j,n}^t|^{pl-1}
\right]^{1/2}.
\end{split}
\end{equation*}
For $l=1$, we take $ p=2$ and
for $l=2,4$ we take $p=3/2$ 
to have
\begin{equation*}
E\left[
\sum_{j=0}^{\infty} |[M]_{j,n}^t|^{pl-1}
\right]^{1/(2p)} E\left[
\sum_{j=0}^{\infty} \left| [M]_{j,n}^t \right|^q
\right]^{1/q} = O(\epsilon_n^l)
\end{equation*}
by the assumption (\ref{add}).
Since $\delta$ can be arbitrarily small, 
this estimate ensures that we can replace 
$\bar{M}$ and $\gamma$ by $\gamma^\delta \cdot M$ and 
$\gamma^\delta$ respectively in
(\ref{Fconv}) and (\ref{Fconv2}).
Put
\begin{equation*}
\gamma^{\delta,n}_s = \gamma^\delta_{\tau^n_j} \text{ for } s \in 
[\tau^n_j, \tau^n_{j+1}).
\end{equation*}
By the same argument, we can estimate
\begin{equation*}
E \left[\epsilon_n^{-l}\sum_{j=0}^{\infty}
\int_{\tau^n_j \wedge t}^{\tau^n_{j+1} \wedge t} \left|
\int_{\tau^n_j \wedge t}^{s}
(\gamma^\delta_u-\gamma^{\delta,n}_u)\mathrm{d}M_u 
\right|^{l} \mathrm{d}\langle M \rangle_s \right]
\end{equation*}
to ensure that  we can replace $(\bar{M} -\bar{M}^n)$ 
with $\gamma^{\delta,n}(M-M^n)$ in
(\ref{Fconv}) and (\ref{Fconv2}), where $M^n$ is defined
as
\begin{equation*}
M^n_s = M_{\tau^n_j},
 \text{ for } j\geq0 \text{ with } s \in [\tau^n_j, \tau^n_{j+1}).
\end{equation*}
Then, repeat the proof for 
the case $\gamma \equiv 1$ with $H$ replaced with 
$H \gamma^\delta$.
\hfill////

{\it Proof of Theorem~\ref{thm1}: } 
Put $\bar{M} = \gamma \cdot M$, 
$A = \psi \cdot \langle M \rangle$ and
define $A^n, \bar{M}^n$ as
(\ref{bars}) with $H = A$.
Then we have
\begin{equation*}
Z^n = (A-A^n)\cdot Y + ((\bar{M}-\bar{M}^n) \varphi) \cdot \langle M \rangle
+ (\bar{M}-\bar{M}^n)\cdot M^Y.
\end{equation*}
We shall  prove that
\begin{equation*}
\epsilon_n^{-1}(A-A^n)\cdot Y_t  \to 0, \ \ 
\epsilon_n^{-1}((\bar{M}-\bar{M}^n) \varphi) \cdot \langle M \rangle_t \to 
\frac{1}{3}(b \varphi) \cdot \langle M \rangle_t
\end{equation*}
in probability uniformly in  $t$ on compact sets of $[0,T)$ and that
\begin{equation} \label{dn0}
D^n := \epsilon_n^{-1}(\bar{M}-\bar{M}^n)\cdot M^Y
\end{equation}
converges $\mathcal{F}$-stably to
\begin{equation*}
\frac{1}{3}(b\gamma) \cdot M^Y + \frac{1}{\sqrt{6}}
(c\gamma) \cdot Y^{\prime}.
\end{equation*}

Step a) Let us show
\begin{equation} \label{stepa}
\epsilon_n^{-1}(A-A^n)\cdot Y_v  \to 0
\end{equation}
uniformly in $v \in [0,t]$ in probability.
Fix $\delta_1, \delta_2>0$  arbitrarily and take
a bounded cag process  
$\psi^\delta$ such that
\begin{equation*}
P[|\psi-\psi^\delta|^2 \cdot \langle M \rangle_t > \delta_1] <  \delta_2
\end{equation*}
by Lemma~\ref{cag}.
Observe that for any $v \in [0,t]$,
\begin{equation} \label{aa}
\begin{split}
&\epsilon_n^{-1}\sum_{j=0}^{\infty}
\int_{\tau^n_j \wedge v}^{\tau^n_{j+1} \wedge v}
\int_{\tau^n_j \wedge v}^{v}\psi_u \mathrm{d}\langle M \rangle_u 
\varphi_s \mathrm{d}\langle M \rangle_s \\
&= \epsilon_n^{-1}\sum_{j=0}^{\infty}
\int_{\tau^n_j \wedge v}^{\tau^n_{j+1} \wedge v}
\int_{\tau^n_j \wedge v}^{s}\psi^\delta_u \mathrm{d}\langle M \rangle_u 
\varphi_s \mathrm{d}\langle M \rangle_s \\
&+ \epsilon_n^{-1}\sum_{j=0}^{\infty}
\int_{\tau^n_j \wedge v}^{\tau^n_{j+1} \wedge v}
\int_{\tau^n_j \wedge v}^{s}(\psi_u -\psi^\delta_u )
\mathrm{d}\langle M \rangle_u 
\varphi_s \mathrm{d}\langle M \rangle_s
\end{split}
\end{equation}
and that 
\begin{equation*}
\begin{split}
&\epsilon_n^{-1}\sum_{j=0}^{\infty}
\int_{\tau^n_j \wedge t}^{\tau^n_{j+1} \wedge t}
\int_{\tau^n_j \wedge t}^{s}|\psi^\delta_u| \mathrm{d}\langle M \rangle_u 
|\varphi_s| \mathrm{d}\langle M \rangle_s \\
& \leq C^\delta \epsilon_n^{-1}\sum_{j=0}^{\infty}
[M]^t_{j,n} \int_{\tau^n_j \wedge t}^{\tau^n_{j+1} \wedge t}
|\varphi_u| \mathrm{d}\langle M \rangle_u \\
&\leq C^\delta \epsilon_n^{-1}\sum_{j=0}^{\infty}
\left\{|[M]^t_{j,n}|^3 \int_{\tau^n_j \wedge t}^{\tau^n_{j+1} \wedge t}
|\varphi_u|^2 \mathrm{d}\langle M \rangle_u\right\}^{1/2} \\
&\leq C^\delta \left\{\epsilon_n^{-2}\sum_{j=0}^{\infty}
|[M]^t_{j,n}|^3
\right\}^{1/2}
\left\{ \int_0^t |\varphi_u|^2 \mathrm{d}\langle M \rangle_u\right\}^{1/2}
\end{split}
\end{equation*}
for a constant $C^\delta$. Using Lemmas~\ref{moreless}, \ref{sup}
and the Burkholder-Davis-Gundy inequality, we have
\begin{equation*}
\epsilon_n^{-2}\sum_{j=0}^{\infty}
|[M]^t_{j,n}|^3 \to 0
\end{equation*}
since
\begin{equation*}
\epsilon_n^{-2}\sum_{j=0}^{N[\tau^n]_t}
G_{j,n}^6 \to 0, \ \ 
\epsilon_n^{-4}\sum_{j=0}^{N[\tau^n]_t}
G_{j,n}^{12} \to 0 
\end{equation*}
in probability by Condition~\ref{gs}. 
Hence the first term of the right hand side of (\ref{aa}) converges to
$0$ uniformly in $v \in [0,t]$ in probability.
For the second term, we have
\begin{equation*}
\begin{split}
&\epsilon_n^{-1}\sum_{j=0}^{\infty}
\int_{\tau^n_j \wedge t}^{\tau^n_{j+1} \wedge t}
\int_{\tau^n_j \wedge t}^{s}|\psi_u -\psi^\delta_u|
\mathrm{d}\langle M \rangle_u 
|\varphi_s| \mathrm{d}\langle M \rangle_s \\
&\leq
\epsilon_n^{-1}\sum_{j=0}^{\infty}
\int_{\tau^n_j \wedge t}^{\tau^n_{j+1} \wedge t}
|\psi_u -\psi^\delta_u |\mathrm{d}\langle M \rangle_u 
\int_{\tau^n_j \wedge t}^{\tau^n_{j+1} \wedge t}
|\varphi_s| \mathrm{d}\langle M \rangle_s \\
&\leq
\sqrt{|\psi - \psi^\delta|^2 \cdot \langle M \rangle_t}
\epsilon_n^{-1}\sum_{j=0}^{\infty}
\sqrt{[M]_{j,n}^t }
\int_{\tau^n_j \wedge t}^{\tau^n_{j+1} \wedge t}
|\varphi_s| \mathrm{d}\langle M \rangle_s \\
&\leq
\sqrt{|\psi - \psi^\delta|^2 \cdot \langle M \rangle_t}
\epsilon_n^{-1}\sum_{j=0}^{\infty} \left\{
|[M]_{j,n}^t|^2
\int_{\tau^n_j \wedge t}^{\tau^n_{j+1} \wedge t}
|\varphi_s|^2 \mathrm{d}\langle M \rangle_s \right\}^{1/2}\\
&\leq
\sqrt{|\psi - \psi^\delta|^2 \cdot \langle M \rangle_t}
\left\{
\sum_{j=0}^{\infty} 
\epsilon_n^{-2}|[M]_{j,n}^t|^2 \right\}^{1/2}
\left\{\int_0 ^t|\varphi_s|^2 \mathrm{d}\langle M \rangle_s \right\}^{1/2}.
\end{split}
\end{equation*}
Using again Lemmas~\ref{moreless}, \ref{sup}
and the Burkholder-Davis-Gundy inequality, we have
\begin{equation*}
\epsilon_n^{-2} \sum_{j=0}^{\infty} 
|[M]_{j,n}^t|^2 = O_p(1).
\end{equation*}
Since $\delta_1, \delta_2$ can be arbitrarily small, 
the left hand sum of (\ref{aa}) converges to $0$ 
uniformly in $v \in [0,t]$ in probability.
Next,   observe that
\begin{equation*}
\langle (A-A^n)\cdot M^Y \rangle
= ((A-A^n)^2 \kappa)\cdot \langle M \rangle,
\end{equation*}
so that it suffices for (\ref{stepa}) to prove
\begin{equation*}
\epsilon_n^{-2}((A-A^n)^2\kappa)
\cdot \langle M \rangle_t \to 0
\end{equation*}
in probability, in light of  the Lenglart inequality.
This follows  from
\begin{equation*}
\begin{split}
&\epsilon_n^{-2} \sum_{j=0}^{\infty} 
\int_{\tau^n_j \wedge t}^{\tau^n_{j+1} \wedge t}
\left\{\int_{\tau^n_j \wedge t}^s \psi_u \mathrm{d} 
\langle M \rangle_u \right\}^2
\kappa_s \mathrm{d}\langle M \rangle_s \\
&\leq 
\epsilon_n^{-2} \sum_{j=0}^{\infty} 
\int_{\tau^n_j \wedge t}^{\tau^n_{j+1} \wedge t}
|\psi_u|^2 \mathrm{d}\langle M \rangle_u |[M]_{j,n}^t|
 \int_{\tau^n_j \wedge t}^{\tau^n_{j+1} \wedge t}
\kappa_s \mathrm{d}\langle M \rangle_s \\
&\leq
\sup_{ j \geq 0} \int_{\tau^n_j \wedge t}^{\tau^n_{j+1} \wedge t}
|\psi_u|^2 \mathrm{d}\langle M \rangle_u
\left\{
\epsilon_n^{-4} \sum_{j=0}^{\infty}| [M]_{j,n}^t|^3
\right\}^{1/2} \left\{
\int_0^t \kappa_s^2 \mathrm{d}\langle M \rangle_s 
\right\}^{1/2},
\end{split}
\end{equation*}
and noticing that
\begin{equation*}
\sup_{ j \geq 0} \int_{\tau^n_j \wedge t}^{\tau^n_{j+1} \wedge t}
|\psi_u|^2 \mathrm{d}\langle M \rangle_u  \to 0, \ \ 
\epsilon_n^{-4}\sum_{j=0}^{\infty}|[M]_{j,n}^t|^3 = o_p(1)
\end{equation*}
by Lemmas~\ref{moreless}, \ref{sup} and Condition~\ref{gs}.

Step b) 
Let us show that
\begin{equation}\label{b0}
\epsilon_n^{-1}((\bar{M}-\bar{M}^n) \varphi) \cdot \langle M \rangle_v \to 
\frac{1}{3}(b \varphi \gamma ) \cdot \langle M \rangle_v
\end{equation}
uniformly in $v \in [0,t]$ in probability. 
Fix $\delta_1, \delta_2 > 0$ arbitrarily and 
take a bounded cag adapted process $\varphi^{\delta}$
such that
\begin{equation*}
P[|\varphi-\varphi^\delta|^2 \cdot \langle M \rangle_t > \delta_1] <  \delta_2
\end{equation*}
by Lemma~\ref{cag}.
Notice that
\begin{equation} \label{ba}
|\epsilon_n^{-1}((\bar{M}-\bar{M}^n) (\varphi-\varphi^\delta)) \cdot 
\langle M \rangle_v|
\leq \sqrt{|\varphi-\varphi^\delta|^2 \cdot \langle M \rangle_t}
\sqrt{\epsilon_n^{-2}(\bar{M}-\bar{M}^n)^2\cdot \langle M \rangle_t}
\end{equation}
and
\begin{equation*}
\epsilon_n^{-2}(\bar{M}-\bar{M}^n)^2\cdot \langle M \rangle_t = O_p(1)
\end{equation*}
by Lemma~\ref{unify}. 
Note also that
\begin{equation} \label{be}
|(b\varphi) 
\cdot \langle M \rangle_v - (b \varphi^\delta)\cdot \langle M \rangle_v|
\leq \sqrt{ b^2 \cdot \langle M \rangle_t }
\sqrt{|\varphi- \varphi^\delta|^2 \cdot \langle M \rangle_t}.
\end{equation}
Since $\delta_1, \delta_2$ can be
arbitrarily small, the estimates (\ref{ba}) and (\ref{be})
ensures that 
we can suppose $\varphi$ is a bounded cag adapted 
process without loss of generality.
Then, putting
\begin{equation*}
\varphi^{n}_s = \varphi_{\tau^n_j} \text{ for }
j \geq 0 \text{ with } s \in [\tau^n_j, \tau^n_{j+1}),
\end{equation*}
we have also that
\begin{equation*}
|\epsilon_n^{-1}((\bar{M}-\bar{M}^n) (\varphi-\varphi^{n})) \cdot \langle M
 \rangle_v|
\leq \sqrt{|\varphi-\varphi^{n}|^2 \cdot \langle M \rangle_t}
\sqrt{\epsilon_n^{-2}(\bar{M}-\bar{M}^n)^2\cdot \langle M \rangle_t}.
\end{equation*}
Note that 
\begin{equation*}
|\varphi-\varphi^{n}|^2 \cdot \langle M \rangle_t \to 0
\end{equation*}
in probability as $n \to \infty$
because $\varphi$ is now assumed to be bounded and left continuous.
Applying Lemma~\ref{unify}, we have (\ref{b0}).

Step c)
Let us study the asymptotic distribution of $D^n$ defined as (\ref{dn0}).
Put
\begin{equation*}
\hat{D}^n = D^n - \frac{1}{3}(b\gamma) \cdot M^Y.
\end{equation*}
In light of Theorem~\ref{jac}, it suffices to
show the following convergences in probability.
\begin{enumerate}
\item
\begin{equation*}
\langle \hat{D}^n, M^Y \rangle_t \to 0,
\end{equation*}
\item
\begin{equation*}
\langle \hat{D}^n \rangle_t \to \frac{1}{6} 
(c^2 \gamma^2 \kappa ) \cdot \langle M \rangle_t, \ \ 
\end{equation*}
\item
\begin{equation*}
\langle \hat{D}^n, \hat{M} \rangle_t \to 0,
\end{equation*}
\end{enumerate}
for all $t \in [0,T)$ and 
for all bounded martingale $\hat{M}$ orthogonal to $M^Y$.
The last one is trivial.
In order to see the first convergence, it suffices to see
\begin{equation*}
\langle D^n, M^Y \rangle_t = 
\epsilon_n^{-1}((\bar{M}-\bar{M}^n)\kappa )\cdot \langle M \rangle_t
\to \frac{1}{3}(b \gamma \kappa) \cdot \langle M \rangle_t
\end{equation*}
in probability. This is shown in the same manner as for (\ref{b0}).
In order to see the second convergence, fix $\delta_1, \delta_2 > 0$
arbitrarily and  take a bounded cag process $\kappa^\delta$
such that
\begin{equation*}
P[|\kappa-\kappa^\delta|^2\cdot \langle M \rangle_t > \delta_1] < \delta_2
\end{equation*}
by Lemma~\ref{cag}.
Notice that
\begin{equation}\label{dn}
\langle D^n \rangle = \epsilon_n^{-2}
((\bar{M}-\bar{M}^n)^2 \kappa^\delta) \cdot \langle M \rangle
+
\epsilon_n^{-2}
((\bar{M}-\bar{M}^n)^2 (\kappa -\kappa^\delta)) \cdot \langle M \rangle 
\end{equation}
and the second term is negligible since
\begin{equation*}
\epsilon_n^{-2}((\bar{M}-\bar{M}^n)^2|\kappa -\kappa^\delta|)
 \cdot \langle M \rangle_t  
\leq \sqrt{\epsilon_n^{-4}(\bar{M}-\bar{M}^n)^4\cdot \langle M \rangle_t }
\sqrt{|\kappa - \kappa^\delta|^2 \cdot \langle M \rangle_t}
\end{equation*}
in light of  Lemma~\ref{unify}
and the fact that
$\delta_1, \delta_2$ can be arbitrarily small.
Furthermore, putting
\begin{equation*}
\kappa^{\delta,n}_s = \kappa^\delta_{\tau^n_j} \text{ for } j \geq 0
\text{ with } s \in [\tau^n_j, \tau^n_{j+1}), 
\end{equation*}
we can replace $\kappa^\delta$ with $\kappa^{\delta, n}$
in the first term of (\ref{dn}) by the same argument.
Then, we have from Lemma~\ref{unify} that
\begin{equation*}
\langle D^n \rangle_t \to \frac{1}{6}(a^2 \gamma^2 \kappa) \cdot 
\langle M \rangle_t
\end{equation*}
in probability. Since
\begin{equation*}
\langle (b\gamma) \cdot M^Y \rangle = (b^2\gamma^2 \kappa )
\cdot \langle M \rangle,
\end{equation*}
it remains only to show
\begin{equation*}
\langle D^n, (b\gamma) \cdot M^Y \rangle_t
\to \frac{1}{3}(b^2\gamma^2 \kappa )
\cdot \langle M \rangle_t
\end{equation*}
in probability. 
Since the left hand side is
\begin{equation*}
\epsilon_n^{-1}(
(\bar{M}-\bar{M}^n)b \gamma \kappa ) \cdot \langle M \rangle_t,
\end{equation*}
the convergence follows from the same argument as for (\ref{b0}).
\hfill////  

{\it Proof of Theorem~\ref{thm5}: }
The convergence (\ref{uconv}) follows from
Condition~\ref{gz} and Lemma~\ref{sup}.
The convergence of $U^n Z^n$ in $D[0,T)$ is a consequence of
the fact that the convergence of $Z^n/\epsilon_n$ is stable.
To show (\ref{boundu2}), we first notice that
\begin{equation*}
G_{j,n}^4/G_{j,n}^2 - \frac{3}{4}|G_{j,n}^3/G_{j,n}^2|^2
\geq |G_{j,n}^2 /G_{j,n}^1|^2 \  \text{ a.s.,}
\end{equation*}
which follows from Lemma~\ref{threefourth}.
In light of Lemma~\ref{sup} and Condition~\ref{gz},
this inequality implies
\begin{equation*}
(Hc^2)\cdot \langle M \rangle \geq  
(H \zeta^{-2})\cdot \langle M \rangle \  \text{ a.s.}
\end{equation*}
for any $H \in \mathcal{P}_M^1$.
Thus we have
\begin{equation*}
\begin{split}
V_t &= \frac{1}{6}(c^2\gamma^2)\cdot \langle Y \rangle_t 
|(|u|\zeta)\cdot \langle M \rangle_t|^2 \\
&\geq
\frac{1}{6}(\zeta^{-2}\gamma^2)\cdot \langle Y \rangle_t 
|(|u|\zeta)\cdot \langle M \rangle_t|^2 
\geq
\frac{1}{6}
\left|
(|u\gamma|^{2/3}\kappa^{1/3})\cdot \langle M \rangle_t\right|^3
\end{split}
\end{equation*}
by H$\ddot{\text{o}}$lder's inequality. \hfill////

{\it Proof of Theorem~\ref{thm6}: }
The convergence (\ref{Nconv}) follows from
Condition~\ref{g2} and Lemma~\ref{sup}.
The convergence of $\sqrt{N[\tau^n]} Z^n$ in $D[0,T)$ is a consequence of
the fact that the convergence of $Z^n/\epsilon_n$ is stable.
To show (\ref{bound2}), we first notice that
\begin{equation*}
G_{j,n}^4/G_{j,n}^2 - |G_{j,n}^3/G_{j,n}^2|^2
\geq G_{j,n}^2 \  \text{ a.s.,}
\end{equation*}
which follows from Lemma~\ref{PearsonLem}.
In light of Lemma~\ref{sup} and Condition~\ref{g2},
this inequality implies
\begin{equation*}
(Hc^2)\cdot \langle M \rangle \geq  
(H q^2)\cdot \langle M \rangle \  \text{ a.s.}
\end{equation*}
for any $H \in \mathcal{P}_M^1$.
Thus we have
\begin{equation*}
\begin{split}
V_t &= \frac{1}{6}(c^2\gamma^2)\cdot \langle Y \rangle_t 
\ q^{-2}\cdot \langle M \rangle_t \\
& \geq
\frac{1}{6}(q^2\gamma^2)\cdot \langle Y \rangle_t 
\ q^{-2}\cdot \langle M \rangle_t 
\geq
\frac{1}{6}
(|\gamma| \sqrt{\kappa})\cdot \langle M \rangle_t
\end{split}
\end{equation*}
by the Cauchy-Schwarz inequality. \hfill////

\section{Effective schemes}\label{effective}
Here we give
effective discretization schemes.
Let $(X,Y,M) \in \mathcal{S}$.
For the sake of brevity, we suppose
$T$ is finite  in this section.
Then, by a localization argument, we can suppose
without loss of generality that
there exists $\delta > 0$ such that 
$\langle M \rangle$ is strictly increasing a.s. on $[T-\delta,T)$.
In fact, we can consider a sequence $M^K$ instead of 
$M$ defined as,
for example,
\begin{equation*}
M^K_t = M_{t \wedge \sigma_K} + \hat{W}_{t} - \hat{W}_{t \wedge
 \sigma_K},
\ \ 
\sigma_K = \inf\{t > 0; \langle M \rangle_t \geq K\} \wedge (T-1/K)
\end{equation*}
with $K \to \infty$, 
where $\hat{W}$ is a Brownian motion defined on an extension of
$(\Omega, \mathcal{F},P)$.
Recall that stable convergence is stable against such a localization
procedure.
Then, for any positive sequence $\epsilon_n$
with $\epsilon_n \to 0$ and for any 
$g \in \mathcal{P}_0$ with $g^{-1} \in \mathcal{P}_0$
the sequence of stopping times $\tau^n$ defined as
\begin{equation} \label{taug}
\tau^n_0 = 0, \ \ \tau^n_{j+1} = 
\inf\left\{ t > \tau^n_j; |M_t - M_{\tau^n_j}| = \epsilon_n g_{\tau^n_j}
\right\} \wedge T
\end{equation}
satisfies Conditions~\ref{gz} and  \ref{g2} with
\begin{equation*}
b_s = 0, \ \ a_s^2 = q^2_s = \zeta_s^{-2} = g_s^2,
\end{equation*}
which follows from  a famous property on the exit times of 
one-dimensional continuous local martingales.
\begin{prop}
Let $(X,Y,M) \in \mathcal{S}$ and
$u \in \mathcal{P}_0$.
The lower bound (\ref{boundu2}) is attained by
$\tau^n$ defined as
(\ref{taug}) with
$g = |u|^{1/3}|\gamma^2
\kappa|^{-1/3}$
if $g, g^{-1} \in \mathcal{P}_0$.
\end{prop}
\begin{prop} \label{prop2}
Let $(X,Y,M) \in \mathcal{S}$.
The lower bound (\ref{bound2}) is attained by 
$\tau^n$ defined as
(\ref{taug}) with $g = |\gamma|^{-1/2}\kappa^{-1/4}$
if $g, g^{-1} \in \mathcal{P}_0$.
\end{prop}
Recall that the lower bound (\ref{bound2}) was derived from 
a combined use of
Lemma~\ref{PearsonLem}
and the Cauchy-Schwarz inequality.
Karandikar~\cite{Karandikar1995} studied a scheme which is defined 
as (\ref{taug}) with
$g=1$ and $X$ instead of $M$ to show the almost sure convergence of $Z^n$. 
In case that $\psi$ appeared in Condition~\ref{YN} 
is locally bounded and $\gamma\equiv1$, 
we can suppose $X=M$ in light of
the Girsanov-Maruyama theorem. Then, we can conclude that
Karandikar's scheme is
 superior to the usual time-equidistant one in that
it yields increments of the integrand which 
attain the equality in Lemma~\ref{PearsonLem}.
It is in fact optimal  if $X=Y$. 

Note that Lemma~\ref{PearsonLem} gives a more precise estimate
\begin{equation*}
c^2 \cdot \langle Y \rangle_t
 = \left(a^2 - \frac{2}{3}b^2\right)
 \cdot \langle Y \rangle_t
\geq \left(\frac{1}{3}b^2 + q^2 \right) \langle Y \rangle_t.
\end{equation*}
The following proposition, for example, is easily shown by this estimate.
\begin{prop} \label{skew}
Let $(X,Y,M) \in \mathcal{S}$, $Z^n$ be defined as (\ref{defz})
and  $\beta, \delta \in \mathcal{P}_0$. 
Denote by  $\mathcal{T}(\beta, \delta)$
the set of sequences of schemes $\tau^n$ which satisfies
Condition~\ref{g2} with $b = \beta$ and 
$q^2 = \delta$. 
Then, for all $t \in [0,T)$,
 $Z^n_t/\epsilon_n$ converges to a mixed normal distribution
with the asymptotic conditional mean
\begin{equation*}
\frac{1}{3} (\beta \gamma) \cdot Y_t
\end{equation*}
and the asymptotic conditional variance $V_t$ satisfying
\begin{equation*}
V_t = \frac{1}{6}
(c^2 \gamma^2) \cdot \langle Y \rangle_t \geq
\frac{1}{6}
\left\{\left(\frac{1}{3} \beta^2 + \delta \right)
\gamma^2\right\}\cdot \langle Y \rangle_t \text{ a.s.}.
\end{equation*}
The equality is attained by 
$\tau^n \in \mathcal{T}_1(M) \cap \mathcal{T}_2(M)$ defined as
\begin{equation*}
\tau^n_{j+1} = \inf\left\{
t > \tau^n_j; M_t -M_{\tau^n_j} \geq
\epsilon_n k_{\tau^n_j}\sqrt{\delta_{\tau^n_j}} \text{ or }
 M_t -M_{\tau^n_j} \leq -\epsilon_n
k_{\tau^n_j}^{-1}\sqrt{\delta_{\tau^n_j}}
\right\}
\end{equation*}
with $\tau^n_0=0$, where
\begin{equation*}
k_s = \frac{\beta_s \delta^{-1/2}_s + \sqrt{\beta_s^2\delta^{-1}_s
 + 4}}{2}.
\end{equation*}
\end{prop}

\section{Conservative delta hedging}
This section treats conservative delta hedging of 
Mykland~\cite{Mykland2000} as an example of
 financial applications.
This framework includes the usual delta hedging for
the Black-Scholes model; even for this classical model,
results presented in this section give a new insight
and a new practical technique for hedging derivatives.
Let $S$ be an asset price process and assume that
\begin{equation*}
dS_t = S_t(\mu_t \mathrm{d}t + \sigma_t \mathrm{d}W_t)
\end{equation*}
is satisfied for predictable processes $\mu$ and $\sigma$
and a standard Brownian motion $W$.
Consider hedging an European contingent claim
$f(S_T)$ for a convex function $f$ of polynomial growth.
Define a function $p$ as
\begin{equation*}
p(S,R,\Sigma) = e^{-R}\int_{\mathbb{R}}f\left(S\exp\{R-\Sigma/2 +
 \sqrt{\Sigma}z\}\right) \phi(z)\mathrm{d}z
\end{equation*}
where $\phi$ is the standard normal density. 
Changing variable, it can be shown that
\begin{equation} \label{pde}
\frac{\partial p}{\partial \Sigma}
= \frac{1}{2}S^2 \frac{\partial^2 p}{\partial S^2}, \ \ 
\frac{\partial p}{\partial R} = 
S\frac{\partial p}{\partial S}- p.
\end{equation}
Put
\begin{equation*}
\eta_K = \inf\{ t> 0; \langle \log(S) \rangle_t \geq K \}
\end{equation*}
 for $K > 0$, $\tilde{V}_t = e^{-rt}V_t$, $\tilde{S}_t = e^{-rt}S_t$ 
and
\begin{equation*}
\Sigma^K_t = K - \langle \log(S) \rangle_{t \wedge \eta_K}, \ \ 
V_t = p\left(S_t,r(T-t),  \Sigma^K_t \right), \ \ 
\pi_t = \frac{\partial p}{\partial S}
\left(S_t,r(T-t), \Sigma^K_t \right),
\end{equation*}
where $r>0$ is a risk-free rate. 
Then, It$\hat{\text{o}}$'s formula and (\ref{pde}) yield 
\begin{equation*}
\tilde{V}_{t \wedge \eta_K} = \int_0^{t \wedge \eta_K}\pi_u 
\mathrm{d}\tilde{S}_u
\end{equation*}
for $t \in [0,T]$, 
that is, the portfolio strategy  $(\pi^0,\pi)$ with
$\pi^0_t = e^{-rt}(V_t- \pi_tS_t)$ is self-financing up to 
$\eta_K \wedge T$.
Moreover, the convexity of $f$ and (\ref{pde})
imply that $p$ is increasing in $\Sigma$, so that
\begin{equation*}
V_T \geq p(S_T,0,0) = f(S_T) \text{ on } \{\eta_K \geq T\}.
\end{equation*}
Note that $p$ is the Black-Scholes price with
cumulative volatility $K$ and that
$\pi$ is the corresponding delta hedging strategy.
The above inequality ensures that the delta hedging super-replicates
any European contingent claim with convex payoff on the set 
$\{\eta_K\geq T\}$.
As $K \to \infty$, $P[\eta_K \geq T] \to 1$, so that
a hedge error due to the incompleteness of
market converges to $0$. Contracts such as variance swap
serve as insurances against the event $\eta_K < T$ for 
predetermined $K$ which is not so large.
See Mykland~\cite{Mykland2003} for an improvement of this conservative
delta hedging.
The purpose here is, however, not to treat such a hedge error due to
the incompleteness but to treat a hedge error due to the restriction
that trades are executed finitely many times in practice.
Note that the rebalancing of a portfolio is  usually 
executed a few times per day while observation of $S$
is almost continuous. Hence, the estimation error of 
$\langle \log(S)\rangle_t$ appeared in $\Sigma_t$ is
negligible compared to the discrete hedging error.
Suppose for brevity that $\eta_K \geq T$ a.s..
An approximation of the strategy  $\pi$ is $\pi^n$ defined as
\begin{equation*}
\pi^n_s = \pi_{\tau^n_j} \text{ for } s \in [\tau^n_j,\tau^n_{j+1})
\end{equation*}
for a discretization scheme $\tau^n$.  
In this context, $N[\tau^n]_t$ is the number of transactions
 up to time $t < T$.
The discounted replication error is given as
\begin{equation*}
Z^n_t = e^{-rt}(V_t- V^n_t) = \int_0^t(\pi_u -\pi^n_u)\mathrm{d}\tilde{S}_u.
\end{equation*}
Notice that after a 
Girsanov-Maruyama transformation,
$X = \pi$ is a local martingale and 
$(X,Y,X) = (\pi,\tilde{S},\pi) \in \mathcal{S}$.
According to  our results in the preceding section, 
the lower bound of the asymptotic variance of 
$\sqrt{N[\tau^n]_t}Z^n_t$ is attained by
the scheme
\begin{equation*}
\tau^n_0 = 0, \ \ 
\tau^n_{j+1} = \inf\left\{
t > \tau^n_j; |\pi_t - \pi_{\tau^n_j}|^2 = \epsilon^2_n 
e^{r\tau^n_j}\Gamma_{\tau^n_j}
\right\} \wedge T,
\end{equation*}
where
\begin{equation*}
\Gamma_t = \frac{\partial^2 p}{\partial S^2}(S_t,r(T-t),\Sigma_t).
\end{equation*}
Note that $\Gamma$ is what is called gamma in financial practice.
In this case, $\tau^n \in \mathcal{T}_1(X)$ with
\begin{equation*}
b_s=0, \ \ 
a_s^2 = c_s^2 = q_s^2 = e^{rs}\Gamma_s,
\end{equation*}
so that we have
\begin{equation*}
 Z \left\{\frac{1}{6}
\int_0^t
e^{-ru} \Gamma_u \mathrm{d} \langle S \rangle_u  \right\}^{1/2}
\end{equation*}
as the asymptotic distribution of $Z^n_t/\epsilon_n$, 
where $Z$ is an independent standard normal variable. 
This scheme is efficient in that the conditional mean-squared
replication error
is asymptotically minimized for a conditionally given number of 
transactions.
The number of transactions $N[\tau^n]_t$ is of course random;
it is high if the path of $\Gamma$ is of high level
because 
$|\pi_t -\pi_{\tau^n_j}|^2/\Gamma_{\tau^n_j} 
\approx |S_t - S_{\tau^n_{j}}|^2\Gamma_{\tau^n_j}$.
This property is intuitively expected  in practice. Note that
$\epsilon_n$ controls the expected number of transactions.
The asymptotic distribution of $\sqrt{N[\tau^n]_t} Z^n_t$ is
\begin{equation*}
\frac{Z}{\sqrt{6}} \int_0^t
e^{-ru} \Gamma_u \mathrm{d} \langle S \rangle_u.
\end{equation*}

In the equidistant case $\tau^n_j = j/n$,
we can apply Theorem~\ref{thm1} to $\epsilon_n = 1/\sqrt{n}$,
$M = W$, $\gamma=\Gamma \sigma S$,
$Y=\tilde{S}$ with
\begin{equation*}
b_s =0, \ \ q_s^2 = 1, \ \ a_s^2 = c_s^2 = 3, \ \ 
N[\tau^n]_t = [nt]
\end{equation*}
to have that  $\sqrt{N[\tau^n]_t}Z^n_t$ 
converges $\mathcal{F}$-stably to
\begin{equation*}
Z \left\{
\frac{t}{2}\int_0^t e^{-2ru}\Gamma_u^2 \sigma_u^2 S_u^2 
\mathrm{d} \langle S \rangle_u
\right\}^{1/2}.
\end{equation*}
The inequality for the asymptotic conditional variance
\begin{equation*}
\frac{1}{6} \left\{
\int_0^t
e^{-ru} \Gamma_u \mathrm{d} \langle S \rangle_u\right\}^2
\leq
\frac{t}{2}\int_0^t e^{-2ru}\Gamma_u^2 \sigma_u^2 S_u^2 
\mathrm{d} \langle S \rangle_u
\end{equation*}
follows directly from the Cauchy-Schwarz inequality.

Karandikar's scheme is defined as
\begin{equation*}
\tau^n_0 = 0, \ \ 
\tau^n_{j+1} = \inf\left\{
t > \tau^n_j; |\pi_t - \pi_{\tau^n_j}| = \epsilon_n 
\right\} \wedge T.
\end{equation*}
After the Girsanov-Maruyama transformation, we 
apply Theorem~\ref{thm1} to $X = M = \pi$ with
\begin{equation*}
b_s = 0, \ \ a_s^2 = c_s^2 = q_s^2 = 1
\end{equation*}
to have that  $\{\sqrt{N[\tau^n]_t}Z^n_t\}$ converges $\mathcal{F}$-stably to
\begin{equation*}
Z \left\{ \int_0^t \Gamma_u^2 \mathrm{d} \langle S \rangle_u\right\}^{1/2}
\left\{
\frac{1}{6}\int_0^t e^{-2ru}\mathrm{d} \langle S \rangle_u
\right\}^{1/2}.
\end{equation*}
The inequality for the asymptotic conditional variance
\begin{equation*}
\frac{1}{6} \left\{
\int_0^t
e^{-ru} \Gamma_u \mathrm{d} \langle S \rangle_u\right\}^2
\leq
 \left\{ \int_0^t \Gamma_u^2 \mathrm{d} \langle S \rangle_u\right\}
\left\{
\frac{1}{6}\int_0^t e^{-2ru}\mathrm{d} \langle S \rangle_u
\right\}
\end{equation*}
follows again directly from the Cauchy-Schwarz inequality.

Taking the purpose of hedging into consideration, 
it might be preferable to use such a scheme $\tau^n$ that
the asymptotic mean of $Z^n_t/\epsilon_n$
is negative. Proposition~\ref{skew} presents
an effective scheme for a given asymptotic conditional mean
and a given asymptotic conditional number of transactions.

More importantly, we can incorporate linear transaction costs.
Suppose that the total cost of
 the delta hedging with a discretization scheme $\tau^n$ up to 
time $t < T$ is 
proportional to
\begin{equation*}
C^n_t = 
\sum_{j=0}^{\infty}|\pi_{\tau^n_{j+1}\wedge t}-\pi_{\tau^n_j \wedge t}|
S_{\tau^n_{j+1}\wedge t}.
\end{equation*}
Let us study the asymptotic distribution of $C^n_t Z^n_t$.
After the Girsanov-Maruyama transformation,
$\pi$ is a local martingale as before.
Notice that
\begin{equation*}
\epsilon_n C^n_t = \epsilon_n \sum_{j=0}^{\infty}
S_{\tau^n_j}
|\pi_{\tau^n_{j+1}\wedge t}-\pi_{\tau^n_j \wedge t}| + o_p(1)
\end{equation*}
if $\tau^n \in \mathcal{T}_1(M)$ with $M = \pi$.
Apply Theorem~\ref{thm5} to $X = M = \pi$, 
$Y = \tilde{S}$
and $u = S$ to have that
$\{C^n_t Z^n_t\}$ converges $\mathcal{F}$-stably to
\begin{equation*}
(S\zeta) \cdot \langle S \rangle_t
\left\{
\frac{1}{3}b  \cdot \tilde{S}_t + Z \sqrt{\frac{1}{6}
c^2 \cdot \langle
\tilde{S} \rangle_t}
\right\} 
\end{equation*}
and that the asymptotic conditional variance of $\{C^n_tZ^n_t\}$
has a lower bound
\begin{equation*}
\frac{1}{6}\left|
\int_0^t|e^{-ru}S_u\Gamma_u^2|^{2/3} \mathrm{d}\langle S \rangle_u
\right|^3,
\end{equation*}
which is attained by $\tau^n$ defined as
\begin{equation*}
\tau^n_0 = 0, \ \ 
\tau^n_{j+1} = \inf\left\{
t > \tau^n_j; |\pi_t -\pi_{\tau^n_j}|^3 = \epsilon_n^3 e^{2r\tau^n_j}
S_{\tau^n_j}\Gamma_{\tau^n_j}^2
\right\} \wedge T.
\end{equation*}
This scheme is efficient in that
the conditional mean-squared replication error
is asymptotically minimized for 
a conditionally given amount of linear transaction costs.

\section{Euler-Maruyama approximation}
Here we propose alternative discretization schemes
for the Euler-Maruyama approximation as an application.
Let us consider the stochastic differential equation
\begin{equation*}
\begin{split}
&\mathrm{d}\Xi_t = \mu(\Xi_t,\eta_t)\mathrm{d}t 
+ \sigma(\Xi_t,\eta_t) \mathrm{d}W_t,\\
&\mathrm{d}\eta_t = \theta(\eta_t)\mathrm{d}t,
\end{split}
\end{equation*}
where $W$ is a one-dimensional 
standard Brownian motion and $\mu, \sigma, \theta$
are continuously differentiable functions.
Since it is rarely possible to generate a path of $\Xi$ fast and exactly,
the Euler-Maruyama scheme is widely used to
approximate to $\Xi$ in simulation.
For sequences $\tau^n = \{ \tau^n_j\}$ with (\ref{tau}), 
the Euler-Maruyama approximation
$\Xi^n$ of $\Xi$ is given as
\begin{equation*}
\begin{split}
&\mathrm{d}\Xi_t^n = \mu(\bar{\Xi}^n_t,\bar{\eta}^n_t)\mathrm{d}t
+ \sigma(\bar{\Xi}^n_t,\bar{\eta}^n_t) \mathrm{d}W_t,\\
&\mathrm{d}\eta_t^n = \theta(\bar{\eta}^n_t)\mathrm{d}t,
\end{split}
\end{equation*}
where $\bar{\Xi}^n_t= \Xi^n_{\tau^n_j}$,
$\bar{\eta}^n_t = \eta^n_{\tau^n_j}$ 
for $j \geq 0$ with $t \in [\tau_j^n, \tau_{j+1}^n)$.
Usually $\tau^n_j = j/n$ is taken.
The convergence rate of the approximation has been extensively
investigated;
see e.g., Kloeden and Platen~\cite{KlPl} for a well-known
strong approximation theorem and
Kohatsu-Higa~\cite{Kohatsu}, Bally and Talay~\cite{BT1,BT2},
Konakov and Mammen~\cite{KM} for
weak approximation theorems.
Newton~\cite{Newton1990} treated passage times.
Cambanis and Hu~\cite{CH} studied
efficiency of deterministic nonequidistant scheme.
Hofmann, M$\ddot{\text{u}}$ller-Gronbach and Ritter~\cite{HMR}
treated a class of adaptive schemes.
Here we exploit a result of
Kurtz and Protter~\cite{KP}, Jacod and Protter~\cite{JP}
to deal with the asymptotic distribution of
pathwise error.  Our aim here is to construct
discretization schemes which are more efficient than the usual
equidistant sampling scheme.

With the aid of localization,
we can suppose $\mu, \sigma$, $1/\sigma$, $\theta$ and their derivatives 
to be bounded.
Suppose that $\tau^n \in \mathcal{T}(W)$.
By Theorem~\ref{thm1}, there exists a conditionally Gaussian
martingale $Z$ such that
\begin{equation*}
\epsilon_n^{-1} (W-W^n) \cdot W \to Z
\end{equation*}
$\mathcal{F}$-stably, where $W^n_t = W_{\tau^n_j}$ for $j \geq 0$ with
$t \in [\tau^n_j, \tau^n_{j+1})$.
Put $L_t^n = \epsilon_n^{-1}(\Xi_t^n - \Xi)$.
Then, applying Kurtz and Protter~\cite{KP},
we have that $L^n$ converges to a process $L$ which satisfies
\begin{equation*}
\mathrm{d}L_t = \partial_1 \mu(\Xi_t,\eta_t)L_t \mathrm{d}t +
\partial_1\sigma(\Xi_t,\eta_t)\left[
L_t\mathrm{d}W_t - \sigma(\Xi_t,\eta_t)\mathrm{d}Z_t
\right],
\end{equation*}
where $\partial_1$ refers to the differential operator with respect to
the first argument.
Solving this stochastic differential equation, we obtain
\begin{equation*}
L_t = -e_t \int_0^t e_s^{-1}\sigma(\Xi_s,\eta_s)
\partial_1\sigma(\Xi_s,\eta_s)\left[
\mathrm{d}Z_s - \partial_1\sigma(\Xi_s,\eta_s) \mathrm{d} 
\langle Z,W \rangle_s
\right],
\end{equation*}
where
\begin{equation*}
e_t = \exp\left\{
\int_0^t \partial_1 \mu(\Xi_s,\eta_s)\mathrm{d}s
+ \int_0^t \partial_1\sigma(\Xi_s,\eta_s)\mathrm{d}W_s
- \frac{1}{2}\int_0^t \partial_1\sigma(\Xi_s,\eta_s)^2\mathrm{d}s
\right\}.
\end{equation*}
Therefore, in light of Theorem~\ref{thm1},
 the distribution of $L_t$ is mixed normal with conditional
mean
\begin{equation*}
-\frac{1}{3}
e_t \int_0^t e_s^{-1}\sigma(\Xi_s,\eta_s)
\partial_1\sigma(\Xi_s,\eta_s) b_s\left[
\mathrm{d}W_t - \partial_1\sigma(\Xi_t,\eta_t) \mathrm{d}t
\right]
\end{equation*}
and conditional variance
\begin{equation} \label{convar}
\frac{1}{6}
e_t^2 \int_0^t e_s^{-2}\sigma(\Xi_s,\eta_s)^2
\partial_1\sigma(\Xi_s,\eta_s)^2 c_s^2\mathrm{d}t.
\end{equation}

\begin{prop}
For any $T > 0$, 
the space-equidistant scheme $\tau^n_{\mathrm{sp}}$ defined as 
(\ref{taug})
with $M=W$, $\epsilon_n = n^{-1/2}$, $g=1$
is three times efficient than the usual time-equidistant scheme
$\tau^n_{\mathrm{tm}} = \{j/n\}$ in the following sense; for any 
$t\in [0,T)$,
\begin{itemize}
\item $E[N^n_t] \leq nt$ and $N^n_t/n \to t$ in probability as
$n\to \infty$ for 
both $N^n = N[\tau^n_{\mathrm{sp}}]$ and
$N^n = N[\tau^n_{\mathrm{tm}}]$,
\item the asymptotic conditional mean  of $L^n_t$ is $0$ for the both schemes,
\item the asymptotic conditional variance of $L^n_t$ for $\tau^n_{\mathrm{sp}}$ is one third of
that for $\tau^n_{\mathrm{tm}}$.
\end{itemize}
\end{prop}
{\it Proof: }
For the space-equidistant case,
\begin{equation*}
b_s \equiv 0, \ \ a_s^2 \equiv c_s^2  \equiv 1, \ \ P[N_t^n]  = n 
P\left[
\sum_{j=0}^{N^n_t-1} |W_{\tau^n_{j+1}} - W_{\tau^n_j}|^2 
\right] 
\leq nt,
\end{equation*}
while
\begin{equation*}
b_s \equiv 0, \ \ a_s^2 \equiv  c_s^2  \equiv 3, \ \ N_t^n = [nt]
\end{equation*}
for the time-equidistant case. \hfill////

Newton~\cite{Newton1990} studied
this  space-equidistant sampling scheme;
the superiority of this scheme is more-or-less known.
The above simple fact of asymptotic conditional variance, however,
has not been recognized so far.
The assumption that $W$ is one-dimensional is a serious restriction.
Nevertheless, also for a stochastic volatility model
\begin{equation*}
\begin{split}
&\mathrm{d}\Xi_t = \hat{\mu}(t,\Xi_t) \mathrm{d}t + \hat{\sigma}(t,V_t) 
[\rho(t,V_t) \mathrm{d}W^1_t + \sqrt{1-\rho(t,V_t)^2}\mathrm{d}W^2_t]
\\& \mathrm{d}V_t = \mu(t,V_t)\mathrm{d}t + \sigma(t,V_t)\mathrm{d}W^1_t
\end{split}
\end{equation*}
with a two-dimensional standard Brownian motion $(W^1,W^2)$,
a scheme defined as (\ref{taug}) with $M=W^1$ 
and $g=1$ results in 
a one third conditional asymptotic variance 
of the  Euler-Maruyama approximation error for $\Xi$.
This is because in light of Theorem~\ref{thm1},
discretization error is determined by only
conditional moments of increments of integrand, which is a
function of $V$ independent of $W^2$ in this example.

Next, let us consider to minimize (\ref{convar})
in case that $\partial_1 \sigma$ is nondegenerate.
Define $\tau^n$ as
\begin{equation}\label{adap}
\tau_0^n =0,  \ \ \tau_{j+1}^n = \inf\{ t > \tau_j^n;
|W_{t} - W_{\tau_j^n}|^2 =
\epsilon(\tau_j^n)
\},
\end{equation}
where
\begin{equation*}
\epsilon(\tau_j^n) = \frac{ \epsilon_n^2 \hat{e}_{\tau_j^n}}
{\sigma(\Xi_{\tau_j^n}^n,
\eta_{\tau_j^n}^n)\partial_1 \sigma(\Xi_{\tau_j^n}^n,
\eta_{\tau_j^n}^n)}
\end{equation*}
and
\begin{equation*}
\begin{split}
\log(\hat{e}_{\tau_j^n}) = \sum_{i=0}^{j-1} & \bigl\{
\partial_1 \mu(\Xi_{\tau_i^n}^n,\eta_{\tau_i^n}^n) (\tau_{i+1}^n -
\tau_i^n)  + \partial_1 \sigma(\Xi_{\tau_i^n}^n,\eta_{\tau_i^n}^n)
(W_{\tau_{i+1}^n} - W_{\tau_i^n}) \\
&- \frac{1}{2}  \partial_1 \sigma(\Xi_{\tau_i^n}^n,\eta_{\tau_i^n}^n)^2
(\tau_{i+1}^n -
\tau_i^n) \bigr\}.
\end{split}
\end{equation*}
Then, Condition~\ref{g2} is satisfied with
\begin{equation*}
b_s \equiv 0, \ \ a_s^2 = c_s^2 = q_s^2 = 
\frac{e_s}{\sigma(\Xi_s,\eta_s)\partial_1\sigma(\Xi_s,\eta_s)}.
\end{equation*}
Proposition~\ref{prop2} implies
that this adaptive scheme attains  a lower bound
for (\ref{convar}) among $\tau^n \in \mathcal{T}(W)$.
In this sense, this scheme is optimal.
A disadvantage of this scheme is the difficulty to
estimate the expected number of data.
In other words, we cannot answer how to
choose $\epsilon_n$ so that the expected number of data is less than
$n$. In practice, it will be better to use
\begin{equation*}
\tau_0^n =0,  \ \ \tau_{j+1}^n = \inf\{ t > \tau_j^n;
|W_{t} - W_{\tau_j^n}|^2 =
\epsilon(\tau_j^n) \vee \epsilon_n^{\prime}
\},
\end{equation*}
for $\epsilon_n^{\prime} > 0$
in order to ensure that  a simulation is done in a finite time.

We conclude this section by a remark on
generating the random variable $(\tau, W_{\tau})$
satisfying
\begin{equation*}
\tau = \inf\{t > 0; |W_t-W_0| = \epsilon\}
\end{equation*}
on a computer for a given $\epsilon$.
It is sufficient that
the distribution function $F_{\epsilon}$ of $\tau$ is available because
\begin{equation*}
P[W_{\tau} = W_0 \pm \epsilon] = 1/2
\end{equation*}
and $\tau, W_{\tau}$ are independent. In fact,
for a random variable $U$ uniformly distributed on $(0,1)$,
\begin{equation*}
(\tau, W_{\tau}-W_0) \sim (F_{\epsilon}^{-1}(2U-[2U]), \epsilon(2[2U]-1)).
\end{equation*}
It is known that the density of $\tau$ is given by
\begin{equation*}
\frac{2}{\sqrt{2\pi t^3}}\sum_{n=-\infty}^{\infty}
(4n+1)\epsilon \exp\left\{-\frac{(4n+1)^2\epsilon^2}{2t}\right\}
\end{equation*}
See Karatzas and Shreve~\cite{KS}, 2.8.11.
Using the fact that
\begin{equation*}
\int_0^t \frac{\alpha}{\sqrt{2\pi t^3}}e^{-\alpha^2/2t} \mathrm{d}t
= 2\int_{\alpha/\sqrt{t}}^{\infty}\phi(x)\mathrm{d}x
\end{equation*}
for $\alpha > 0$, we obtain $F_{\epsilon}(t) = G(\epsilon/\sqrt{t})$, where
\begin{equation*}
\label{gfunc}
G(x) = 4 \sum_{n=0}^{\infty}
(\Phi((4n+3)x)- \Phi((4n+1)x) ).
\end{equation*}
According to our numerical study,
$G(x) \approx 1$ for $0 \leq x \leq 0.1$.
This is not surprising because
$G(0+) = 1$, $G^{\prime}(0+) = G^{\prime \prime}(0+)=0$.
Note that If $x \geq 0.1$,
the speed of convergence of the infinite series is very fast.
We can therefore use
\begin{equation*}
G(x) \approx \begin{cases}
4 \sum_{n=0}^{\lfloor N/x \rfloor}(
\Phi((4n+3)x)- \Phi((4n+1)x) )
& x \geq 0.1, \\
1  & 0 \leq x < 0.1
\end{cases}
\end{equation*}
for, say,  $N=3$ as a valid approximation of $G$.
It is noteworthy
that $G$ is independent of $\epsilon$, so that
once we obtain the inverse function of $G$ numerically,
it is done very fast to generate $\tau$ repeatedly even if
$\epsilon$ changes adaptively as in (\ref{adap}).
Note also that 
\begin{equation*}
G(x) \leq 4(1-\Phi(x)), \ \
G^{-1}(y) \leq \Phi^{-1}(1-y/4).
\end{equation*}
These inequalities will be useful in
numerical calculation of $G^{-1}$ for sufficiently
small $y$ (~large $x$~). Besides, if $x \geq 3$,
$G(x) \approx 4(1-\Phi(x))$ and
$G^{-1}(y) \approx \Phi^{-1}(1-y/4)$.


\begin{appendix}

\section{Auxiliary results}
Here we give auxiliary results for the proof of our theorems.
The following limit theorem, which plays an essential role in
this article,  is a 
simplified version of a result of Jacod~\cite{Jacod1994} and
Jacod and Shiryaev~\cite{JS}, Theorem  IX.7.3,
which extends a result of  Rootz\'en~\cite{Rootzen}.
Let $M = \{ M_t, \mathcal{F}_t, 0 \leq t < \infty  \}$ be a 
continuous local martingale defined on $(\Omega, \mathcal{F}, P)$
and  $\mathcal{M}^\perp$ be the set of
bounded 
$\{\mathcal{F}_t\}$-martingales orthogonal to $M$.

\begin{thm}\label{jac} 
Let $\{Z^n\}$ be a sequence of continuous $\{\mathcal{F}_t\}$-local
 martingales. Suppose that  there exist
an $\{\mathcal{F}_t\}$-adapted continuous process 
$V = \{V_t\}$ such that  for all  $\hat{M} \in \mathcal{M}^\perp $, 
$ t \in [0,\infty)$, 
\\
\begin{equation*}
\langle Z^n , \hat{M} \rangle_t \to 0,
\ \ 
\langle Z^n, M \rangle_t \to 0, \ \ 
\langle Z^n \rangle_t \to V_t
\end{equation*}
in probability. 
Then, the $C[0,\infty)$-valued sequence $\{Z^n\}$ converges
$\mathcal{F}$-stably to 
the distribution of the time-changed process $W^{\prime}_V$
where $W^{\prime}$ is a standard Brownian motion independent of
$\mathcal{F}$.
\end{thm}

The following lemma is repeatedly used in our proofs.
\begin{lem} \label{moreless}
Consider a sequence of filtrations
\begin{equation*}
\mathcal{H}_j^n \subset \mathcal{H}_{j+1}^n, \ \ 
j,n \in \mathbb{Z}_+ = \{0,1,2,\dots\}
\end{equation*}
and random variables $\{ U_j^n \}_{j \in \mathbb{N}}$ 
with  $U_j^n$ being  $\mathcal{H}_j^n$-measurable.
Let $N^n(\lambda)$ be a $\{\mathcal{H}_j^n\}$-stopping time for 
each $n \in \mathbb{Z}_+$ and $\lambda$ which
is an element of a set $ \Lambda$.
Let  $U(\lambda)$ be a random variable for each $\lambda \in \Lambda$.
If it holds that there exists $\lambda_0 \in \Lambda$ such that
\begin{equation*}
N^n(\lambda) \leq N^n(\lambda_0) \text{ a.s. for all } \lambda \in \Lambda
\text{ and }
\sum_{j=1}^{N^n(\lambda_0)} P[ |U_j^n|^2 | \mathcal{H}_{j-1}^n ] 
\to 0 
\end{equation*}
as $n \to \infty$, 
then the following two are equivalent;
\begin{enumerate}
\item
\begin{equation*}
\sup_{\lambda \in \Lambda} \left|\sum_{j=1}^{N^n(\lambda)} U_j^n -
 U(\lambda) \right| \to 0 \  \text{ as } n \to \infty.
\end{equation*}
\item
\begin{equation*}
\sup_{\lambda \in \Lambda} \left|\sum_{j=1}^{N^n(\lambda)} 
P[U_j^n | \mathcal{H}_{j-1}^n ] -  U(\lambda) \right|  \to 0
\ \text{ as } n \to \infty.
\end{equation*}
\end{enumerate}
Here the convergences
are in probability.
\end{lem}
{\it Proof: }
Note that
\begin{equation*}
\sup_{\lambda \in \Lambda}
\left|
\sum_{j=1}^{N^n(\lambda)}U^n_j - 
\sum_{j=1}^{N^n(\lambda)}E[U^n_j|\mathcal{H}_{j-1}^n]
\right|
\leq \sup_{k \in \mathbb{N}}
\left| \sum_{j=1}^{k} V_j^n\right|,
\end{equation*}
where
\begin{equation*}
V_j^n = (U_j^n - E[U^n_j|\mathcal{H}_{j-1}^n])1_{j\leq N^n(\lambda_0)}.
\end{equation*}
By the Lenglart inequality, we have
\begin{equation*}
P\left[
\sup_{k \in \mathbb{N}}
\left| \sum_{j=1}^{k} V_j^n\right| \geq \epsilon
\right] \leq 
\frac{\eta}{\epsilon^2} + P\left[
\sum_{j=1}^{\infty}E[|V_j^n|^2|\mathcal{H}_{j-1}^n] \geq \eta
\right]
\end{equation*}
for any $\epsilon, \eta > 0$.
The result then follows from the convergence
\begin{equation*}
\sum_{j=1}^{\infty}E[|V_j^n|^2|\mathcal{H}_{j-1}^n]
\leq \sum_{j=1}^{N^n(\lambda_0)}E[|U^n_j|^2 | \mathcal{H}_{j-1}^n]
\to 0
\end{equation*}
in probability.
\hfill ////

The following lemma is well-known, so its proof is omitted.

\begin{lem}\label{cag}
For all  $H \in \mathcal{P}_M^k$, $t \in [0,T)$ 
and $\delta_1, \delta_2 > 0$, there
exists a bounded adapted left-continuous process $\hat{H}$ such that 
\begin{equation*}
P[|H - \hat{H}|^k \cdot \langle M \rangle_t > \delta_1] <  \delta_2.
\end{equation*}
\end{lem}

The following lemma is taken from Fukasawa~\cite{F2009}.
\begin{lem} \label{sup}
Let $M$ be a continuous local martingale
with $E[\langle M \rangle_T^6]  < \infty$ and suppose that
$\tau^n \in \mathcal{T}(M)$. Then, for all  $t\in [0,T) $,
\begin{equation} \label{sup1}
\sup_{j \geq 0, 0 \leq s \leq t} 
|M_{\tau_{j+1}^n \wedge s} -M_{ \tau_j^n \wedge s}|^2 
= o_p(\epsilon_n)
\end{equation}
as well as 
\begin{equation} \label{sup2}
\sup_{0 \leq j \leq N[\tau^n]_t } 
|\langle M \rangle_{\tau_{j+1}^n} - 
\langle M \rangle_{\tau_j^n} | = o_p(\epsilon_n),
\ \ 
\sup_{0 \leq j \leq N[\tau^n]_t} |M_{\tau_{j+1}^n } -M_{ \tau_j^n }|^2 
= o_p(\epsilon_n),
\end{equation}
where $N[\tau^n]_t$ is defined as (\ref{ntau}). In particular,
\begin{equation}\label{ninf}
N[\tau^n]_t 
\to \infty \text{ a.s. }
\end{equation}
as $n \to \infty$. Moreover, for all locally bounded adapted cag process $f$, it holds
\begin{equation}\label{fconv}
\sum_{j=0}^{N[\tau^n]_t} f_{\tau^n_j} G_{j,n}^2 \to 
\int_0^t f_s \mathrm{d}\langle M \rangle_s
\end{equation}
in probability, uniformly in $t$ on compact sets of $[0,T)$.
\end{lem}
{\it Proof: } 
Note that
\begin{equation*}
\sum_{j=0}^{N[\tau^n]_t} G_{j,n}^2 = O_p(1)
\end{equation*}
since $E[\langle M \rangle_T] < \infty$.
By the assumptions, it follows 
\begin{equation*}
\sum_{j=0}^{N[\tau^n]_t} E[(M_{\tau_{j+1}^n} - M_{\tau_j^n})^{2k}| 
\mathcal{F}_{\tau_j^n}] = o_p(\epsilon_n^k)
\end{equation*}
for $k=3,6$ and with the aid of Lemma~\ref{moreless}, we have
\begin{equation*}
\sum_{j=0}^{N[\tau^n]_t} (M_{\tau_{j+1}^n} - M_{\tau_j^n})^6 
= o_p(\epsilon_n^3).
\end{equation*}
On the other hand,
\begin{equation*}
\sup_{0 \leq j \leq N[\tau^n]_t}|M_{\tau_{j+1}^n}- M_{\tau_j^n}|^2 
\leq \left\{
\sum_{j=0}^{N[\tau^n]_t} (M_{\tau_{j+1}^n} - M_{\tau_j^n})^6 
\right\}^{1/3},
\end{equation*}
so that the second of (\ref{sup2}) follows.

To show 
(\ref{sup1}), 
we use Doob's maximal inequality to have
\begin{equation*}
E\left[ \sup_{0 \leq t <  \infty}
|M_{\tau_{j+1}^n \wedge t} - M_{\tau_j^n \wedge t}|^{2k}|
\mathcal{F}_{\tau_j^n }\right] / G_{j,n}^2
= o_p(\epsilon_n^k)
\end{equation*}
for $k=3,6$. Using Lemma~\ref{moreless} again, we obtain
\begin{equation*}
\sum_{j=0}^{N[\tau^n]_t}
\sup_{0 \leq s <  \infty}
|M_{\tau_{j+1}^n \wedge s} - M_{\tau_j^n \wedge s}|^{6}
= o_p(\epsilon_n^3),
\end{equation*}
which implies 
(\ref{sup1}) since
\begin{equation*}
\sup_{j \geq 0, 0\leq s \leq t} 
|M_{\tau_{j+1}^n \wedge s} -M_{ \tau_j^n \wedge s}|^2 
\leq \left\{\sum_{j=0}^{N[\tau^n]_t}
\sup_{0 \leq s <  \infty}
|M_{\tau_{j+1}^n \wedge s} - M_{\tau_j^n \wedge s}|^{6} \right\}^{1/3}.
\end{equation*}
Using the Burkholder-Davis-Gundy inequality and
Doob's maximal inequality, we have also
\begin{equation*}
\sum_{j=0}^{N[\tau^n]_t} E[|\langle M \rangle_{\tau_{j+1}^n} - 
\langle M \rangle_{\tau_j^n}|^k|
\mathcal{F}_{\tau_j^n}] = o_p(\epsilon_n^k)
\end{equation*}
for $k=3,6$, which implies the first of (\ref{sup2}) in the same manner.
Note that (\ref{ninf}) follows from
\begin{equation*}
N[\tau^n]_t  + 1
\geq \frac{\langle M \rangle_t}
{\sup_{0 \leq j \leq N[\tau^n]_t } 
|\langle M \rangle_{\tau_{j+1}^n} - 
\langle M \rangle_{\tau_j^n} |}.
\end{equation*}
To see (\ref{fconv}), again in light of Lemma~\ref{moreless},
it suffices to observe that
\begin{equation*}
\sum_{j=0}^{N[\tau^n]_t} f_{\tau^n_j}(\langle M \rangle_{\tau^n_{j+1}} - 
\langle M \rangle_{\tau^n_{j}}) \to 
\int_0^t f_s \mathrm{d}\langle M \rangle_s
\end{equation*}
and
\begin{equation*}
\sum_{j=0}^{N[\tau^n]_t} f_{\tau^n_j}^2 G_{j,n}^4  = 
\epsilon_n^2 
\sum_{j=0}^{N[\tau^n]_t} f_{\tau^n_j}^2 a_{\tau^n_j}^2G_{j,n}^2 + 
o_p(\epsilon_n^2)
\to 0.
\end{equation*}
\hfill////

\section{Kurtosis-Skewness inequality}
\begin{lem} \label{PearsonLem}
Let $X$ be a random variable with $E[X] = 0$ and $E[X^4] < \infty$.
Then it holds
\begin{equation} \label{pearineq}
\frac{E[X^4]}{E[X^2]^2} - \frac{E[X^3]^2}{E[X^2]^3} \geq 1.
\end{equation}
The equality is attained if and only if $X$ is a Bernoulli
random variable.
\end{lem}
{\it Proof: }
This is called Pearson's inequality and shown easily as follows.
\begin{equation*}
E[X^3]^2 = 
E[X(X^2-E[X^2])]^2 \leq E[X^2]E[|X^2-E[X^2]|^2]
= E[X^2](E[X^4]-E[X^2]^2).
\end{equation*}
If the equality holds, then $X$ and $X^2-E[X^2]$ must be linearly
dependent, so that $X$ must be Bernoulli.
\hfill////

The following lemma gives a similar inequality to the above.
The proof is however rather different and the result itself
is seemingly new.

\begin{lem} \label{threefourth}
Let $X$ be a random variable with $E[X] = 0$ and $E[X^4] < \infty$.
Then it holds
\begin{equation} \label{thfo}
\frac{E[X^4]}{E[X^2]^2} - \frac{3}{4}\frac{E[X^3]^2}{E[X^2]^3} \geq 
\frac{E[X^2]}{E[|X|]^2}.
\end{equation}
The equality is attained if and only if $X$ is a Bernoulli
random variable.
\end{lem}
{\it Proof: } We divide the proof into 4 steps.

{\bf Step a) } 
It is straightforward to see that
the equality holds if
$X$ is a Bernoulli random variable with $E[X]=0$.

{\bf Step b) } Let us show if $E[X] =0$ and 
the support of $X$ is a finite set,
then the distribution $P^X$ of $X$  is a finite mixture of
Bernoulli distributions with mean $0$.
First, consider the case $n=3$.
Suppose without loss of generality that
\begin{equation*}
P[X=a] = p, \ 
P[X=b] = q, \ 
P[X=c] = r, \ \ p+q+r=1
\end{equation*}
with $a > b \geq 0 > c$. Put
\begin{equation*}
P_1(a) = \frac{-c}{a-c}, \ \ P_1(c) = \frac{a}{a-c}, \ \ 
P_2(b) = \frac{-c}{b-c}, \ \ P_2(c) = \frac{b}{b-c}.
\end{equation*}
Then $P_1$ and $P_2$ define Bernoulli distributions with mean $0$
and supports $\{a, c\}$ and $\{b, c\}$ respectively.
Putting $\lambda = (c-a)p/c = P^X(a)/P_1(a)$, we have
\begin{equation*}
\lambda P_1(a) = p, \ \ 
(1-\lambda) P_2(b) = q, \ \ 
\lambda P_1(c) +  (1-\lambda) P_2(c) = r,
\end{equation*}
which means
\begin{equation*}
P^X = \lambda P_1 +  (1-\lambda) P_2.
\end{equation*}
Now, let us treat the general case by induction.
Suppose that the claim holds for the case of $n$
and consider the case of $n+1$. Without loss of generality,
we suppose 
\begin{equation*}
P[X = a_j] = p_j, \ \ j=0, 1, \dots, n, \ \ 
p_0 + p_1 + \cdots + p_n = 1
\end{equation*}
with
\begin{equation*}
a_0 > a_1 > \cdots > a_k \geq 0 > a_{k+1} > \cdots > a_n
\end{equation*}
for an integer $k$, $2 \leq k \leq n-1$. Put
\begin{equation*}
\tilde{a}_1 = \frac{a_0p_0 + a_1p_1}{p_0 + p_1}, \ \ 
\tilde{p}_1 = p_0 + p_1, \ \ 
\tilde{a}_j = a_j, \ \ 
\tilde{p}_j = p_j, \ \ 2 \leq j \leq n
\end{equation*}
and 
\begin{equation*}
\tilde{P}(\tilde{a}_j) = \tilde{p}_j, \ j=1, \dots n.
\end{equation*}
Notice that $\tilde{P}$ defines a distribution with mean $0$ 
and support $\{\tilde{a}_1,a_2, \dots, a_n\}$. By the assumption of
induction, there exists $\tilde{\lambda}_{ij} \geq 0$, 
$1 \leq i \leq k$, $k < j \leq n$ such that
\begin{equation*}
\sum_{i,j} \tilde{\lambda}_{ij} = 1, \ \
\tilde{P} = \sum_{i,j} \tilde{\lambda}_{ij} \tilde{P}_{ij}, \ \ 
\tilde{P}_{ij}(\tilde{a}_i) = 
\frac{-\tilde{a}_j}{\tilde{a}_i-\tilde{a}_j}, \ \ 
\tilde{P}_{ij}(\tilde{a}_j) =\frac{\tilde{a}_i}{\tilde{a}_i-\tilde{a}_j}.
\end{equation*}
Here $\tilde{P}_{ij}$ defines a Bernoulli distribution
with mean $0$ and support $\{\tilde{a}_i, \tilde{a}_j\}$.
Now consider a distribution $Q_j$ defined as
\begin{equation*}
Q_j(a_0) = \frac{p_0}{p_0+p_1} \tilde{P}_{1j}(\tilde{a}_1), \ \ 
Q_j(a_1) = \frac{p_1}{p_0+p_1} \tilde{P}_{1j}(\tilde{a}_1), \ \ 
Q_j(a_j) = \tilde{P}_{1j}(\tilde{a}_j)
\end{equation*}
for $k < j \leq n$. Notice that $Q_j$ is a distribution with mean $0$
and support $\{a_0, a_1, a_j\}$.
As seen above for the case $n=3$, 
putting  $\mu_j = Q_j(a_0)/P_{0j}(a_0)$,
we have
\begin{equation*}
Q_j = \mu_j P_{0j} + (1- \mu_j) P_{1j},
\end{equation*}
where we define
\begin{equation*}
P_{ij}(a_i) = 
\frac{-a_j}{a_i-a_j}, \ \ 
P_{ij}(a_j)
=\frac{a_i}{a_i-a_j}, \ \ 0 \leq i \leq k, \ k < j \leq n.
\end{equation*}
Putting 
\begin{equation*}
\lambda_{0j} = \mu_j \tilde{\lambda}_{1j}, \ \ 
\lambda_{1j} = (1-\mu_j) \tilde{\lambda}_{1j}, \ \ 
\lambda_{ij} = \tilde{\lambda}_{ij}, \ 2 \leq i \leq k, \ k<j \leq n,
\end{equation*}
we have
\begin{equation*}
\sum_{i,j} \lambda_{ij} = 1 ,\ \ 
P^X = \sum_{i,j} \lambda_{ij} P_{ij},
\end{equation*}
which completes the induction.

{\bf Step c) } Let us show that the function $f(u,v,w,y)$ defined as
\begin{equation} \label{fconc}
f(u,v,w,y) = u - \frac{3}{4}v^2/w - w^3/y^2
\end{equation}
is a concave function. Note that the inequality (\ref{thfo}) follows from
Steps a, b and c since every distribution can be approximated
arbitrarily close by a distribution supported by a finite set.
By a straightforward calculation, the Hessian matrix of $f$ is 
given by
\begin{equation} \label{Hes}
H = \left(
\begin{array}{@{\,}cccc@{\,}}
0 & 0& 0 & 0 \\
0 & -\frac{3}{2w} & \frac{3v}{2w^2} & 0 \\
0 & \frac{3v}{2w^2} & -\frac{3v^2}{2w^3}-\frac{6w}{y^2} &
 \frac{6w^2}{y^3} \\
0 & 0 & \frac{6w^2}{y^3} &  -\frac{6w^3}{y^4}
\end{array}
\right).
\end{equation}
Again by a straightforward calculation, the determinant of 
$H - x I$ is of form $x^2(x+\alpha)(x+\beta)$ with $\alpha > 0$,
$\beta > 0$, which means that $H$ is negative semi-definite.

{\bf Step d) } It remains to show that the equality holds only if
$X$ is a Bernoulli random variable. 
Suppose that there exists a random variable $X$ with $E[X]=0$
such that the equality holds in (\ref{thfo})
which is not Bernoulli.
Recall that the equality holds in (\ref{pearineq})
only if $X$ is Bernoulli. It implies that
the vector of the first four moments of $X$ does not coincide
with that of a Bernoulli random variable.
Note that there exists a random variable $\hat{X}$ of which
the support is a finite set, such that
\begin{equation*}
E[\hat{X}] =E[X], 
\ E[|\hat{X}|] = E[|X|],
\ E[\hat{X}^2] = E[X^2],
\ E[\hat{X}^3] = E[X^3],
\ E[\hat{X}^4] = E[X^4].
\end{equation*}
This can be proved by the Hahn-Banach theorem.
Hence,
we assume the support of $X$ is a finite set
without loss of
generality.
Then, by Steps b and c, there exist Bernoulli distributions 
$P_1$ and $P_2$ and $\lambda \in (0,1)$
such that $P_1 \neq P_2$ and
\begin{equation} \label{f0}
f(\lambda m_1 + (1-\lambda)m_2 ) = 0 
= \lambda f(m_1) + (1 -\lambda)f(m_2),
\end{equation}
where $f$ is defined as (\ref{fconc}) and
\begin{equation*} 
m_i = \left(
\int a^4P_i(da), \int a^3P_i(da), \int a^2P_i(da), \int |a|P_i(da)
\right)^{\prime}, \ \ i=1,2.
\end{equation*}
Here $\prime$ means the transpose of matrix.
By (\ref{Hes}), the eigenvectors of
the Hessian matrix $H$ associated to the eigenvalue $0$ are
\begin{equation*}
h_1 = (1,0,0,0)^\prime, \ \ h_2 = (0,v,w,y)^\prime.
\end{equation*}
Therefore, (\ref{f0}) implies that there exists a constant $c$ 
such that $\bar{m}_2 = c \bar{m}_1$, where
\begin{equation*} 
\bar{m}_i = \left(
\int a^3P_i(da), \int a^2P_i(da), \int |a|P_i(da)
\right)^{\prime}, \ \ i=1,2.
\end{equation*}
It suffices then to show that $\bar{m}_2 = c \bar{m}_1$ 
implies $c = 1$ and that $\bar{m}_1$ uniquely determines
a Bernoulli distribution.
Set
\begin{equation*}
P_2(a) = p, \ P_2(-b) = q, \ p+q=1, \ ap = bq, \ \bar{m}_1 = (v,w,y)^{\prime}
\end{equation*}
and 
\begin{equation*}
a^3p - b^3q = cv, \ \ a^2p + b^2q = cw, \ \  ap + bq = cy.
\end{equation*}
Then we obtain that
\begin{equation*}
a = \frac{2wq}{y}, \ \
b = \frac{2wp}{y},
 \end{equation*}
so that
\begin{equation*}
\frac{2v}{y} = a^2-b^2 = \frac{4w^2}{y^2}(1-2p).
\end{equation*}
Therefore, $a, b, p, q$ are uniquely determined independently of $c$,
which completes the proof.
\hfill//// 
\end{appendix}

\begin{flushleft}
{\bf Acknowledgement: }
This work is partially supported by
JSPS, Grant-in-Aid for Young Scientists (B) and
JST, CREST.
\end{flushleft}

\end{document}